# A Minimalist Approach to Rolling Wheels

Antonín Slavík, Faculty of Mathematics and Physics, Charles University, Prague, Czech Republic ⟨slavik@karlin.mff.cuni.cz⟩

Stan Wagon, Macalester College, St. Paul, MN, USA ⟨wagon@macalester.edu⟩

**Abstract.** In 1960, G. B. Robison discovered the general equations relating roads and wheels, where either can have an unusual shape (e.g., the square wheel rolls smoothly on a catenary). But he used some inobvious assumptions regarding the meaning of *rolling*. Here we derive the equations for the road appropriate for a given wheel using only the single assumption that rolling occurs with no slipping. We do not require that the wheel be differentiable, so this allows the construction of a wheel–road pair when the wheel is a continuous, nowhere differentiable function.

## 1. INTRODUCTION.

The appropriate road for a round wheel is a horizontal line: the wheel rolls smoothly along the line. The converse question is classic: What is the appropriate road for a wheel that is a horizontal line? The answer is a catenary; this is the basis of the famous construction of a road of linked catenaries on which a square wheel will roll smoothly (Figures 1 and 4; a version where the road loops around in a circle is at the National Museum of Mathematics in New York City [**5**]; a recently completed bridge near London, England, is a large square that rolls on giant catenaries [**3**]).

The square wheel problem was posed and solved in 1960 by G. B. Robison [**1**]; see also [**2**]. He discovered general equations relating wheels and roads, but his work is problematic because it is based on certain assumptions as to what *rolling* means. In particular, Robison assumed an "equilibrium property": the wheel's center stays directly above the contact point of the wheel with the road during rolling. He also assumed that the arc lengths of the wheel and road to the point of contact match. These properties are also assumed in [**2**]. One can justify the equilibrium property, but it is not obvious that it should be taken as an axiom. Our goal here is to show how the theory of rolling wheels can be derived solely from the assertion that rolling occurs with no slipping. Thus a teacher who wishes to present the mathematics underlying the square wheel can do so from a single intuitive assumption. By avoiding arc length, our approach does not require that the wheel be differentiable, and so one can imagine a smoothly rolling wheel that arises from a continuous, nowhere differentiable function (Example 7). In fact, the arc length condition by itself is insufficient to guarantee proper rolling; see section 4.

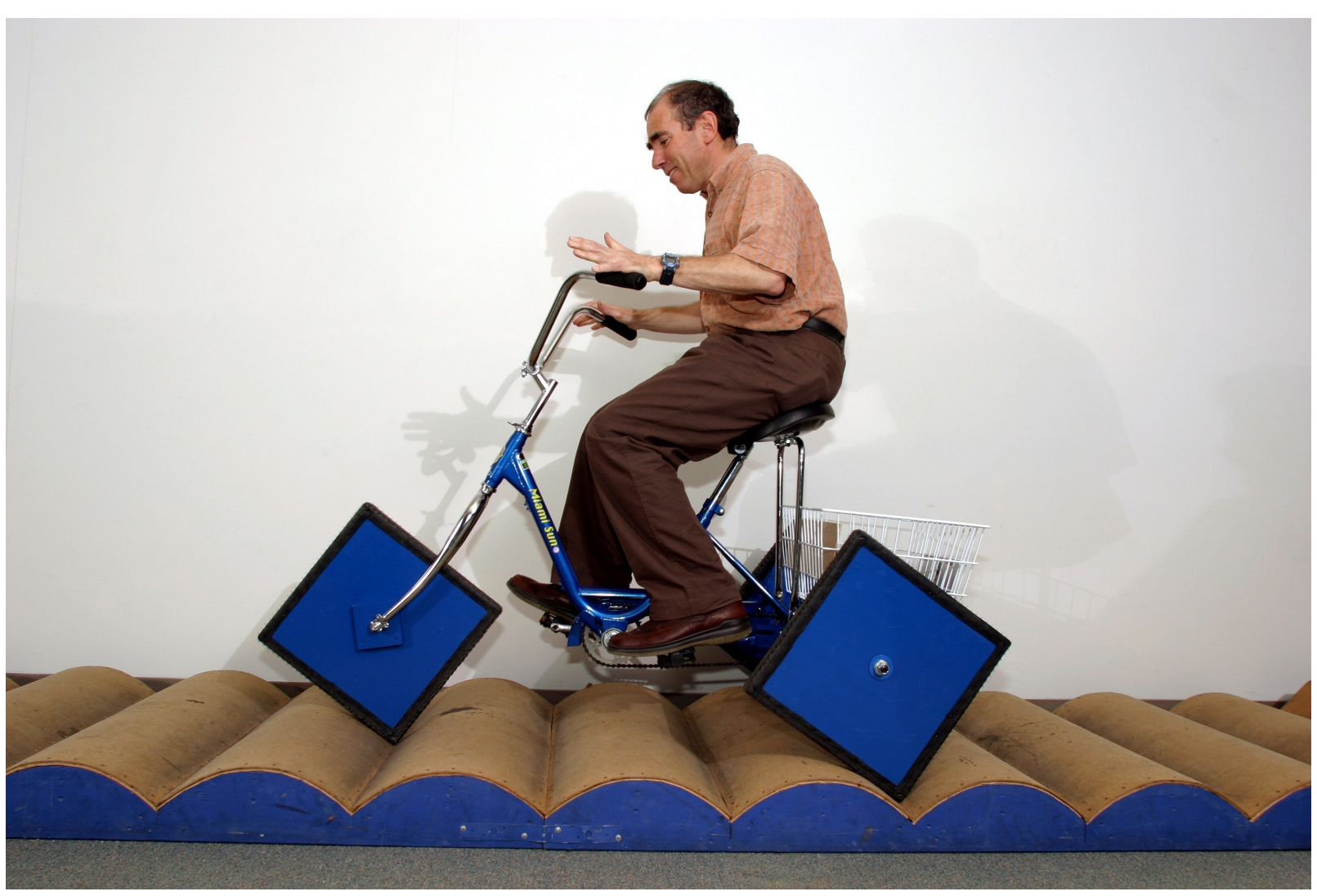

**Figure 1.** Stan Wagon riding the square-wheel tricycle at Macalester College. The ride is smooth because the centers of the wheels, and hence the seat, move only horizontally.

## 2. THE ROAD–WHEEL RELATIONSHIP: $x' = r = -y$.

Suppose we have a wheel $W$, given in polar coordinates in the nonstandard form $W(\theta) = r(\theta)\,(\sin\theta, -\cos\theta)$, with $(0, 0)$ as the specified center. The wheel need not be a closed curve: $r(\theta)$ might not be periodic. Assume that $r$ is continuous and defined on a nondegenerate, possibly infinite, interval $I$ containing 0, and $r(\theta) > 0$ for each $\theta$. Then $W(0) = (0, -r(0))$ is directly below the origin; this is taken to be the initial contact point with the road.

We wish to consider wheels that roll smoothly along a road in that the center of the wheel moves neither up nor down. If the wheel is used on a vehicle such as a bicycle or cart, we do not want the vehicle to bounce up and down. The center can be the geometrical center of the wheel, but it can also be any point that rides along with the wheel. A wheel rolling smoothly along a road from left to right undergoes clockwise rotation around its center followed by horizontal translation as the contact point moves along the road. The horizontal condition enforces the smoothness property. If the center is the physical center of gravity of the wheel, then the smooth rolling means that no work will be done against gravity during the rolling. This can be important: the 27,000-pound square-wheel bridge at Cody Dock in East London can be rolled along giant catenaries, by hand, into an upside-down position so that a boat can pass beneath it [**3**].

Suppose that the wheel has rotated by an angle $\phi$ and let $\tau(\phi)$ be the horizontal displacement associated with this rotation. Then if $P$ is any point on the wheel, its position during rolling is $\sigma_P(\phi) = M(\phi)\cdot P + (\tau(\phi), 0)$, where $M(\phi) = \begin{pmatrix} \cos\phi & \sin\phi \\ -\sin\phi & \cos\phi \end{pmatrix}$ is the matrix for clockwise rotation by $\phi$. Now let $\rho(\theta)$ denote the rotation angle required for a point $P = W(\theta)$ to come into contact with the road (Figure 2). After the corresponding rotation and translation, the position of $P$ becomes $\sigma_P(\rho(\theta)) = M(\rho(\theta))\cdot P + (\tau(\rho(\theta)), 0)$. To summarize, rolling requires the existence of two functions, $\rho(\theta)$ and $\tau(\phi)$, where the domain of $\rho$ is the interval $I$, $\rho(0) = 0$, the domain of $\tau$ is the codomain of $\rho$, $\rho$ is continuous, and $\tau$ is differentiable. We then impose only the following single condition, which will suffice to derive the shape of the road. The main point is that the condition implies that $\rho(\theta) = \theta$.

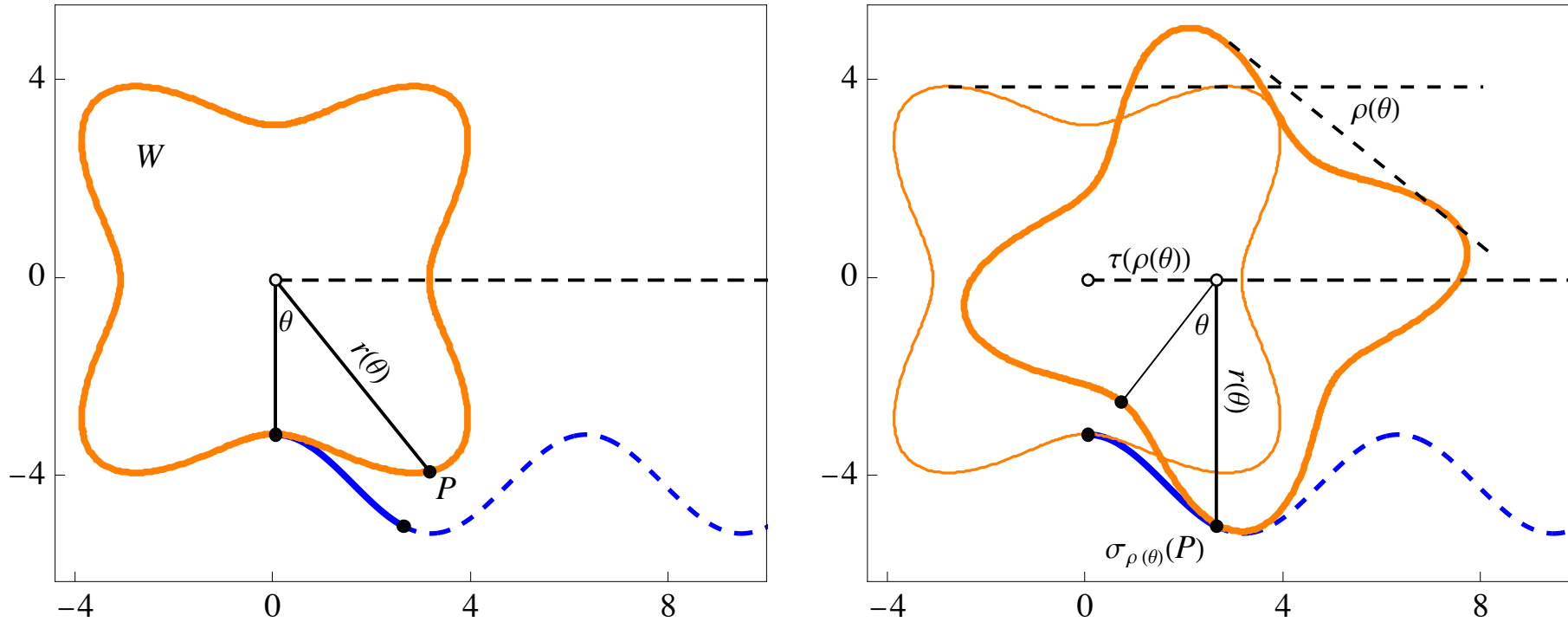

**Figure 2.** The general setup for a rolling wheel given by $r(\theta)$. A consequence of the theorem is that $\rho(\theta)$ equals $\theta$.

**The No-slipping Condition.** During rolling at constant angular speed, each point on the wheel instantaneously stops when it is a contact point. That is, for each $P = W(\theta)$, $\sigma'_P(\rho(\theta)) = (0, 0)$.

The following result gives the location of the contact point, which defines the appropriate road for the given wheel. The concise form of the formula is nice, but not really surprising. The road must conform to the shape of the wheel and, if we accept the aforementioned equilibrium property, that explains $y = -r$. The horizontal distance traveled by a rolling object is proportional to the radius—a larger bicycle wheel moves faster than a smaller one for the same amount of pedalling—and that explains $x' = r$.

**Theorem.** *When a wheel W rolls smoothly along a road, the contact point of* $W(\theta)$ *as the wheel rolls is* $\left(\int_0^\theta r(t)\, dt,\ -r(\theta)\right)$.

*Proof.* We use these simple facts: $M'(\phi) = M\left(\phi + \frac{\pi}{2}\right)$; $M\left(\phi + \frac{\pi}{2}\right)\cdot(0, b) = M(\phi)\cdot(b, 0)$; $M(\alpha)\cdot M(\beta) = M(\alpha+\beta)$; and $W(\theta) = M(-\theta)\,(0, -r(\theta))$. By the no-slipping condition,

$$
\begin{aligned}
(0, 0) &= \sigma'_{W(\theta)}(\rho(\theta)) = M'(\rho(\theta))\cdot W(\theta) + (\tau'(\rho(\theta)), 0) \\
&= M\left(\rho(\theta) + \frac{\pi}{2}\right)\cdot M(-\theta)\cdot(0, -r(\theta)) + (\tau'(\rho(\theta)), 0) = M(-\theta)\cdot M\left(\rho(\theta) + \frac{\pi}{2}\right)\cdot(0, -r(\theta)) + (\tau'(\rho(\theta)), 0) = M(\rho(\theta) - \theta)\cdot(-r(\theta), 0) + (\tau'(\rho(\theta)), 0).
\end{aligned}
$$

To preserve the 0 in the $y$-coordinate, the rotation in $M(\rho(\theta) - \theta)$ must be through either 0 or $\pi$; therefore $\rho(\theta) = \theta + n\pi$. But $\rho$ is continuous and $\rho(0) = 0$, so $\rho(\theta) = \theta$. This means that the contact point is $M(\theta)\cdot W(\theta) + (\tau(\theta), 0) = M(\theta)\cdot M(-\theta)\cdot(0, -r(\theta)) + (\tau(\theta), 0) = (\tau(\theta), -r(\theta))$, which proves the $y$-coordinate of the theorem. Because $\rho(\theta) = \theta$, the $x$-component of the vanishing derivative is $\tau'(\theta) - r(\theta)$ so $\tau'(\theta) = r(\theta)$, which, because the contact point's $x$-coordinate is $\tau(\theta)$, yields the desired integral. □

The theorem implies that the road must have the form $\left(\int_0^\theta r(t)\, dt,\ -r(\theta)\right)$. It is easy to see that this parametrized curve works. Just set $\rho(\theta) = \theta$ and $\tau(\phi) = \int_0^\phi r(t)\, dt$; then

$$
\sigma'_{W(\theta)}(\rho(\theta)) = \sigma'_{W(\theta)}(\theta) = M(\theta - \theta)\cdot(-r(\theta), 0) + (\tau'(\theta), 0) = (-r(\theta), 0) + (r(\theta), 0) = (0, 0).
$$

From now on, we use $(x(\theta), y(\theta))$ for the road for a polar wheel given by $r(\theta)$; then $(x(\theta), y(\theta)) = \left(\int_0^\theta r(t)\, dt,\ -r(\theta)\right)$. A more useful description is that $x(0) = 0$, $x'(\theta) = r(\theta)$, and $y(\theta) = -r(\theta)$. This formulation allows $x(\theta)$ to be obtained by solving a differential equation, often numerically; see Example 7 below. Moreover, $x(\theta)$ is differentiable and strictly increasing, which means that its inverse function $\theta(x)$ exists, is differentiable, and satisfies $\theta'(x) = 1/r(\theta)$.

**Corollary.** ***(a)*** *The road on which a polar wheel rolls smoothly has the functional form* $y = -r(\theta(x))$, *where* $\theta'(x) = -1/r(\theta)$ *and* $\theta(0) = 0$.

***(b)*** *When a wheel smoothly rolls along a road, the center stays directly above the contact point.*

***(c)*** *If the wheel is piecewise continuously differentiable, then so is the road and the arc lengths from the initial contact point to any other contact point are equal.*

***(d)*** *If the wheel is piecewise differentiable, then so is the road and the slopes of each at any point of contact where the slopes exist are equal.*

*Proof.* (a) This is immediate because $\theta(x)$ exists and is monotonic and differentiable, $x'(\theta) = r(\theta)$, and $y(\theta) = -r(\theta)$.

(b) The center after rolling is at $(\tau(\theta), 0)$, which, by the proof, equals $(x(\theta), 0)$, directly above the contact point $(x(\theta), y(\theta))$.

(c) Work on the intervals corresponding to the piecewise property. Because $(x'(\theta), y(\theta)) = (r(\theta), -r(\theta))$, the road is piecewise continuously differentiable. Within each such interval we have $x'(\theta) = r(\theta) = -y(\theta)$, which means $x'(\theta)^2 + y'(\theta)^2 = r(\theta)^2 + r'(\theta)^2$ and the integrals of the square roots of these continuous functions give the arc lengths.

(d) Differentiability is as in (c). Working at the points of differentiability and using the theorem, $\frac{dy}{dx} = \frac{dy/d\theta}{dx/d\theta} = -\frac{r'(\theta)}{r(\theta}$. And, using $\rho(\theta) = \theta$, the slope of the appropriately rotated wheel is the quotient of the terms in $M(\theta) \cdot W'(\theta) = (-r(\theta), r'(\theta))$. This proves equality of slopes. □

The literature [**1**, **2**] derives the road equations by assuming various other conditions (the arc lengths up to the contact point match, the slopes at a contact point match, the center stays directly above the contact point). The corollary shows that these conditions all follow from the no-slipping condition. However, because $r$ need not be differentiable, the arc length of the road might not exist (as in Example 7 below). We have this open question: If we assume only that the wheel is rectifiable, is it true that the road is rectifiable and the arc lengths to the contact point match? Note that the no-slipping condition is equivalent to the assertion that the instantaneous rolling motion is always a clockwise rotation around the point of contact.

The theorem also solves the inverse problem. Given a road $y = f(x)$ whose graph is below the $x$-axis, find the wheel (with center at the origin) that will roll smoothly along it. By part (a) of the corollary, we seek a positive function $r$ such that if we define $x(\theta) = \int_0^\theta r(t)\, dt$ and let $\theta(x)$ be its inverse function (which exists because $x$ is increasing), then $-r(\theta(x)) = f(x)$ for all $x$. Replacing $x$ by $x(\theta)$, we get $f(x(\theta)) = -r(\theta) = -x'(\theta)$. The differential equation $x'(\theta) = -f(x(\theta))$ together with the initial condition $x(0) = 0$ can be solved numerically or symbolically. Once we have $x(\theta)$, we get the wheel in polar form from $r(\theta) = -f(x(\theta))$. The corresponding road is then $-r(\theta(x)) = f(x(\theta(x))) = f(x)$ as required. See Examples 3 and 4.

## 3. EXAMPLES.

Here we present several examples of the method that yields a road for a wheel, or vice versa. Many more can be found in [**1**, **2**].

### 1. Road for a straight line and a square.

Consider the wheel that is the straight line $y = -1$, with center at $(0, 0)$; because $\sec\theta(\sin\theta, -\cos\theta) = (\tan\theta, -1)$, the polar form of the line is $r(\theta) = \sec\theta$, with $-\pi/2 < \theta < \pi/2$. Then, using $u = \sin\theta$, the theorem gives

$$x = \int \sec\theta\, d\theta = \int \frac{\cos\theta}{1-\sin^2\theta}\, d\theta = \int \frac{1}{1-u^2}\, du = \tanh^{-1} u; \text{ and}$$

$$y = -\sec\theta = -\sec\left(\sin^{-1} u\right) = -\frac{1}{\sqrt{1-u^2}} = -\frac{1}{\sqrt{1-\tanh^2 x}} = -\cosh x.$$

Thus the road is an inverted catenary; see Figure 3. Remarkably, the connection between the line and catenary was first discovered by James Clerk Maxwell, who included it in his 1849 paper [**4**], published when he was 18 years old.

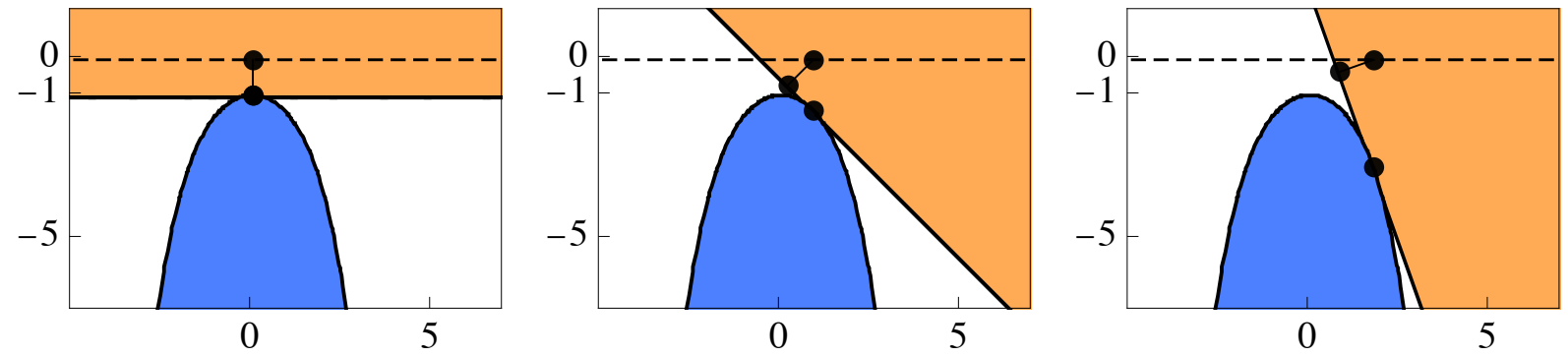

**Figure 3.** An infinite line rolling on the inverted catenary $y = -\cosh x$.

Now we can truncate the line to $|x| \le 1$ and rotate it to form a square centered at the origin. If we then truncate the catenary to $|x| \le \sinh^{-1} 1$ (each slope angle is then 45° from vertical at the ends) and form an infinite track by translating the arches, we will have a road on which a square can roll smoothly (Figure 4).

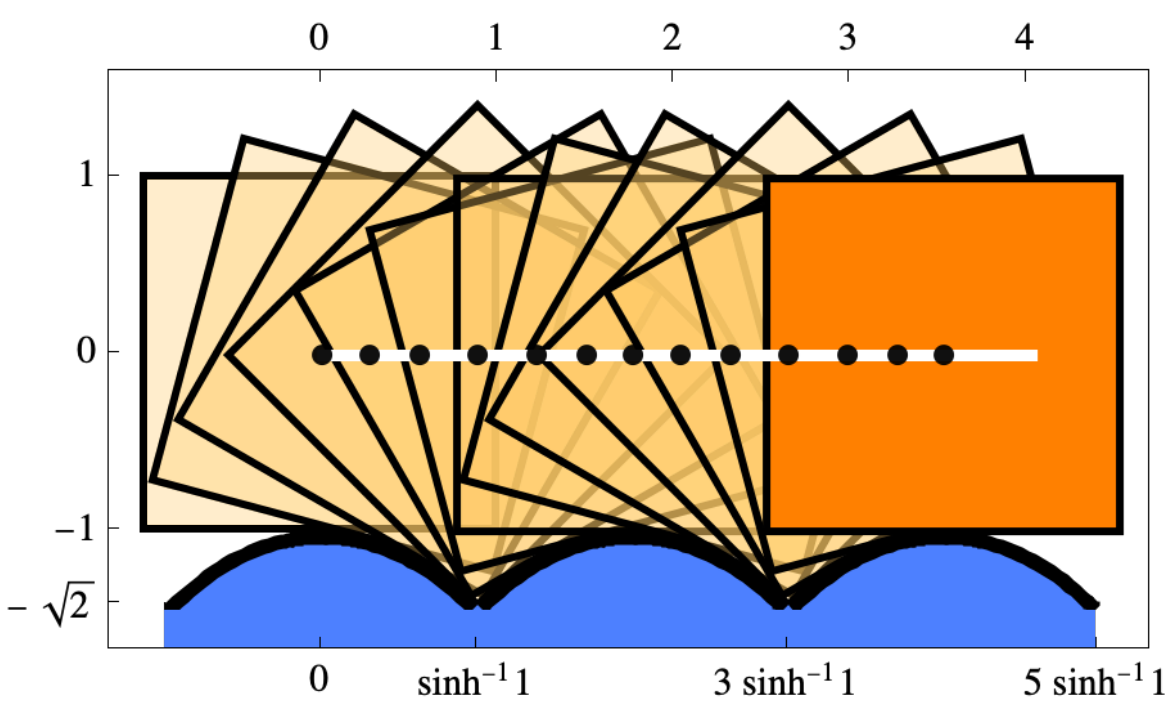

**Figure 4.** A square rolls smoothly on a road made up of truncated inverted catenaries.

As one pedals a square-wheel vehicle at constant angular speed the forward motion is not constant. More precisely, the connection is given by the nonlinear function $x = \cosh^{-1}(\sec \theta)$. Figure 5 shows how $x$ varies with $\theta$.

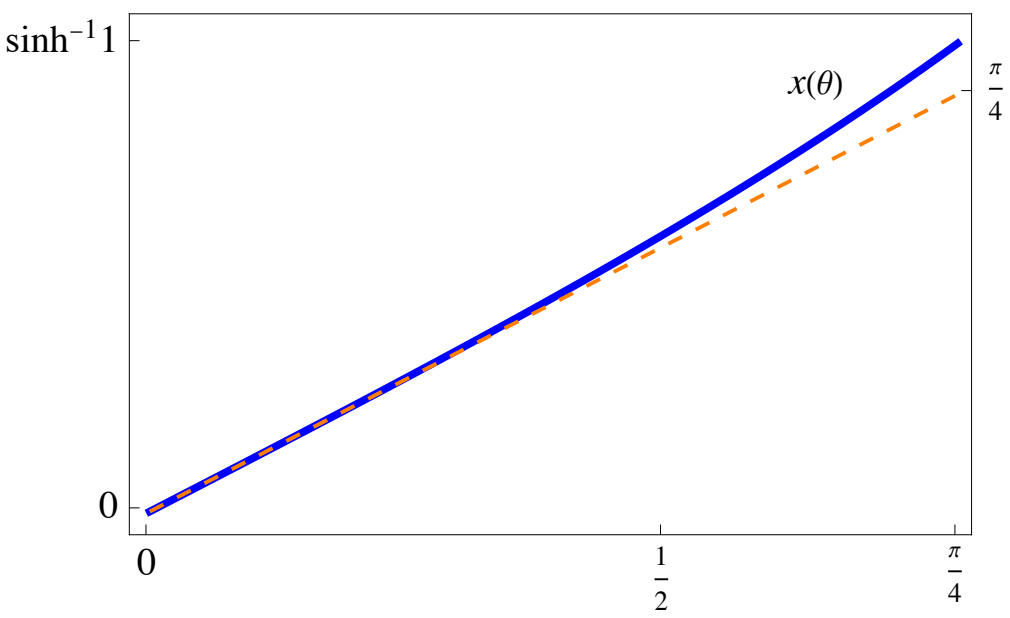

**Figure 5.** The plot of $x$ against $\theta$ for a rolling square.

The preceding idea extends to any wheel that is a regular polygon, but there is one important distinction that arises. A triangular wheel will crash into the road before it reaches the cusp (Figure 6). From a purely geometrical viewpoint there is nothing wrong and a similar situation where collisions are simply ignored arises in Example 7: a continuous but nowhere differentiable wheel.

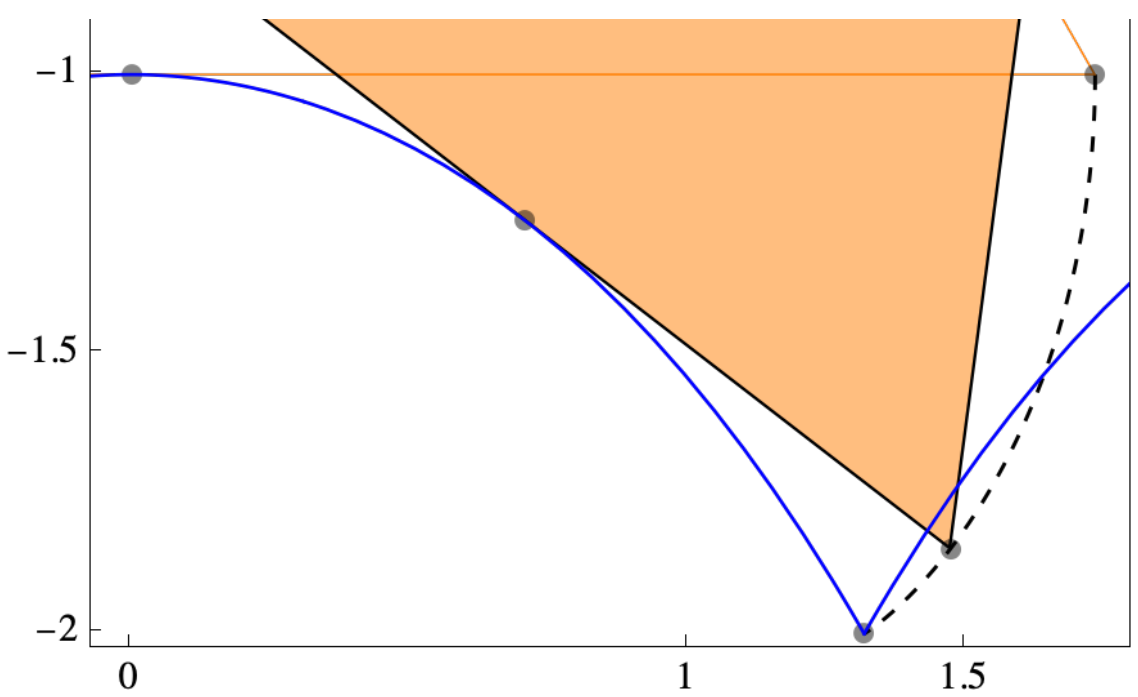

**Figure 6.** A rolling triangle crashes into the road before it reaches the cusp. The motion of a vertex during the rolling (dashed) is the involute of the road.

The many and varied models provide physical proof that the square can roll freely, but a mathematical proof should be provided as well, giving certainty to the fact that a collision as in the triangular case does not happen for a square.

**Proposition.** A square rolls on linked catenaries with no collision.

*Proof*. Let $M_\phi$ and $\tau$ be as above. Let $V$ be the vertex of the square that is initially at $(-1, 1)$ and restrict consideration to the bottom edge of the square and to $0 \leq \phi < \pi/4$; let $P$ be the cusp at $\left(\sinh^{-1} 1, -\sqrt{2}\right)$. The locus of $V$ during the rolling is $\Lambda = M_\phi \cdot V + (\tau(\phi), 0) = \left(\cos\phi - \sin\phi + \cosh^{-1}(\sec\phi), -\sin\phi - \cos\phi\right)$. The second arch of the road stays under the tangent line to this arch $P$ (see Figure 7). But the locus of $V$ is above this line because the dot product $(\Lambda - P) \cdot (1, 1) = \left(\sqrt{2} - 2\sin\phi\right) + \left(\cosh^{-1}(\sec\phi) - \sinh^{-1} 1\right)$ is positive; the first two summands add to a positive number, as do the last two.

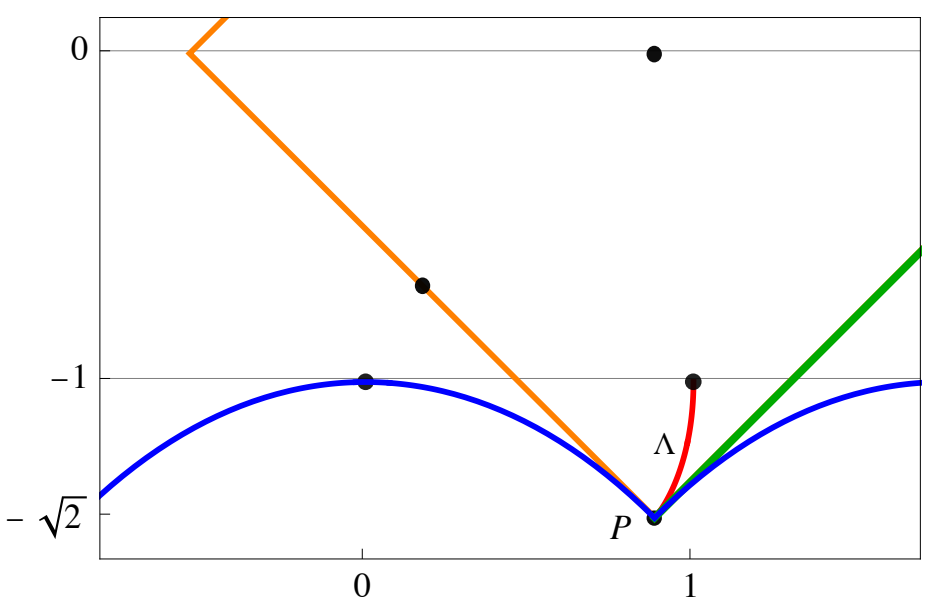

**Figure 7.** The locus of a vertex of a rolling square (red) stays above the tangent line to the second arch of the road, and so a collision cannot occur.

One must also show that no point on the square other than vertex $V$ can crash into the road. Let $Q$ be a point on the bottom edge of the square. Then $Q$ cannot collide with the first arch because it lies on a tangent line to the arch (by Corollary (d)). And $Q$ cannot collide with the second arch because $V$ does not collide and $V$ would have to crash through first. (Let $L$ be the line through $V$ perpendicular to $VQ$. Move $L$ parallel to itself to $L'$, tangent to the arch; the arch is entirely below $L'$ while $LP$ is above it.) It follows from symmetry that the other edges never collide with the road. □

We next prove that the right angle is the critical angle: an acute angle leads to a collision and an obtuse angle does not.

**Proposition.** The vertex of an acute angle rolling on linked catenaries passes below the road in a neighborhood of the cusp. A rolling obtuse angle does not collide with the road.

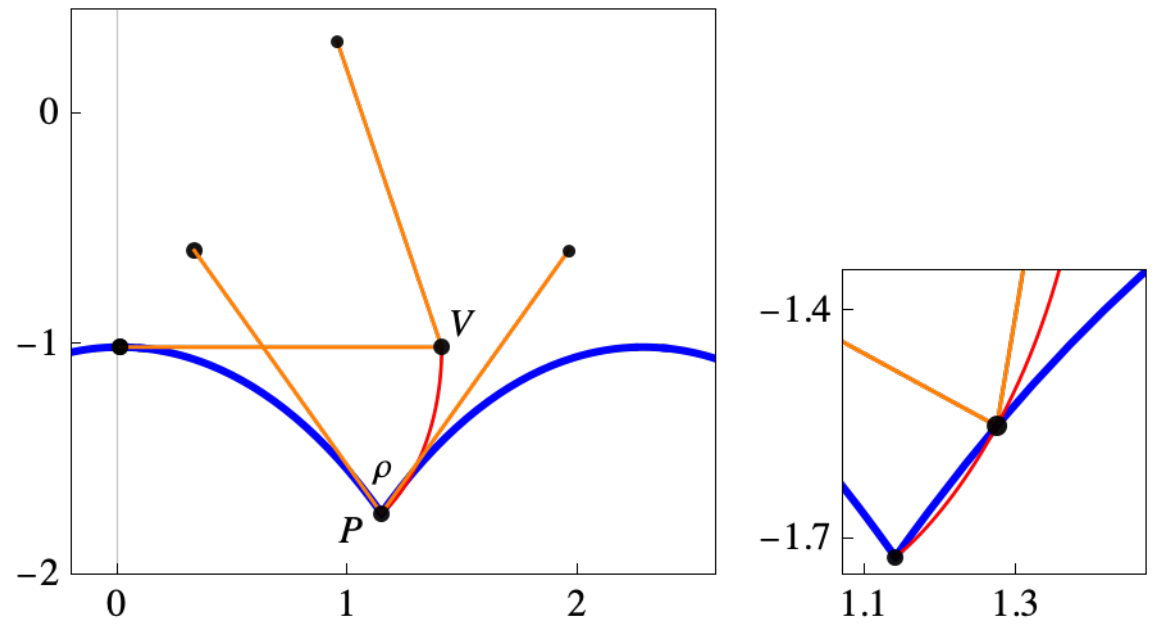

**Figure 7.** A rolling acute angle crashes through the road.

*Proof*. Consider a wheel that is an acute angle $\rho$ with vertex $V$. We may assume that the two legs have equal length and that implies the road consists of two congruent catenaries linked at a cusp $P$ (Figure 7). The angle at $P$ must be $\rho$, for if it were smaller the angle would not fit and if it were larger there would be rotation around $V$ which would cause the center to stray from horizontal. Let $\Lambda$ be the locus of $V$ as the angle rolls into $P$; then $\Lambda$ is the involute of the first catenary from anchor point $(0, -1)$. Therefore we can use the classical result that $\Lambda$ makes a right angle with the catenary at $P$. But $\rho < \pi/2$, and therefore an interval of $\Lambda$ near $P$ lies below the second catenary, implying the desired collision. Alternatively, one can just compute the slope of $\Lambda$ directly from $M_\phi \cdot V + (\tau(\phi), 0)$; it is $\cot\phi$, which gives $\tan(\rho/2)$ for the rotation angle $(\pi - \rho)/2$. This is less than the slope of the second catenary at $P$, which is $\cot(\rho/2)$. These arguments also apply to the obtuse case, showing no collision. □

## A self-contained derivation

It is easy to specialize the theorem to the central case of a line rolling on a catenary, which might be valuable to a teacher wanting to explain only how a square wheel works. Here is a self-contained derivation. Suppose a straight line wheel $W(\theta)$ has polar form $r = \sec\theta$; assume $0 \le \theta < \pi/2$ as symmetry handles the negative case. Let $\rho$ and $\tau$ be as in the definition of rolling: $\rho(\theta)$ gives the clockwise rotation angle for $W(\theta)$ to contact the road and $\tau(\phi)$ is the amount of horizontal translation when the wheel rotates through $\phi$ while rolling. Assume that any point on the line comes to an instantaneous stop when it contacts the road. This means that

$$M'(\rho(\theta)) \cdot W(\theta) + (\tau'(\rho(\theta)), 0) = (\tau'(\rho(\theta)) - \cos(\rho(\theta)) - \sin(\rho(\theta)) \tan\theta, \sin(\rho(\theta) - \theta) \sec\theta) = (0, 0).$$

The vanishing of the $y$-coordinate implies $\rho(\theta) - \theta = n\pi$, and then $\theta = \rho(\theta)$ because $\rho$ is continuous and $\rho(0) = 0$. This means the contact point is $M(\theta) \cdot W(\theta) + (\tau(\theta), 0) = M(\theta) \cdot M(-\theta) \cdot (0, -\sec\theta) + (\tau(\theta), 0) = (\tau(\theta), -\sec\theta)$. The $x$-coordinate of the condition must be 0 so, using $\rho(\theta) = \theta$, one gets $\tau'(\theta) = \sec\theta$, which implies $x = \int \sec\theta \, d\theta$. Integration then yields the road $y = -\cosh x$ as follows:

$$x = \int \sec\theta\, d\theta = \int \frac{\cos\theta}{1-\sin^2\theta}\, d\theta == \int \frac{1}{1-u^2}\, du = \frac{1}{2}\int \frac{1}{1+u} - \frac{1}{u-1}\, du = \tanh^{-1}(\sin\theta) = \cosh^{-1}(\sec\theta) = \cosh^{-1}(-y).$$

## 2. An infinite loop related to the rolling line.

A very simple self-contained approach to the rolling line is to start with the line and catenary and verify that they work. Assume $0 \le \theta < \pi/2$; symmetry yields the negative case. Take $y = -\cosh x = -r(\theta)$ and $r(\theta) = \sec\theta$; then $x = \cosh^{-1}(\sec\theta)$. It is then easy to verify that there is no slipping when $\rho$ and $\tau$ are properly specified: $\rho(\theta) = \theta$ and $\tau(\phi) = \cosh^{-1}(\sec\phi)$. This works because

$$M(\theta)\cdot W(\theta) + (\tau(\theta), 0) = (\tau(\theta), -\cosh(\tau(\theta)),$$

which shows that there is contact when $\rho(\theta) = \theta$, and $M'(\theta)\cdot W(\theta) + (\tau'(\theta), 0) = (0, 0)$, which shows that there is no slipping on contact.

Examining the preceding argument, the crux is that, at the contact point, rotational velocity exactly cancels translational velocity in the $x$-component (it is easy to see that both are 0 in the $y$-component). The rotational speed is $-\sec\theta$ and the translational speed is $\frac{d}{d\theta}\cosh^{-1}(\sec\theta)$, so the key property is that $f(\theta) = \sec\theta$ and $g(\theta) = \cosh^{-1}\phi$ satisfy $\frac{d}{d\theta} g(f(\theta)) = f(\theta)$. Among the trig functions there is one other case of this phenomenon: $\frac{d}{d\theta}\cos^{-1}(\operatorname{sech}\theta) = \operatorname{sech}\theta$. This means that the wheel $r = \operatorname{sech}\theta$ $(-\infty < \theta < \infty)$, a Poinsot spiral, will roll smoothly on the graph of $\cos x$, restricted to $(-\pi/2, \pi/2)$ (Figure 8). Here $x(\theta) = \tan^{-1}(\sinh\theta)$. The $x$-speed approaches 0 as $\theta \to \infty$, so the wheel's center never travels beyond $x = \pi/2$.

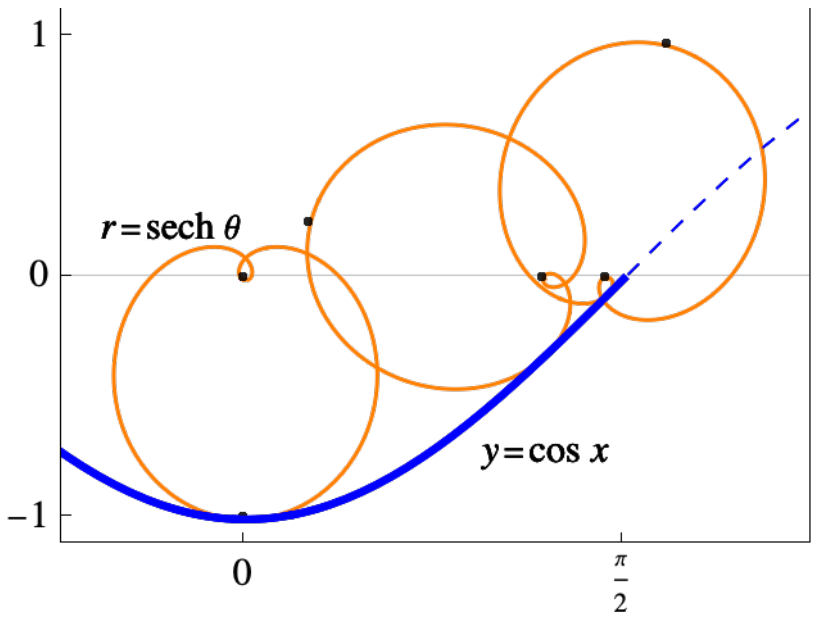

**Figure 8.** A Poinsot spiral winds around the origin infinitely many times. It rolls along a single arch of the cosine function. It is a hyperbolic companion to the line–catenary relationship.

## 3. The inverse square wheel.

Suppose we are given a road defined by the catenary $y = -\cosh x$. Using rotation angle $\theta$, consider this as $(x(\theta), -\cosh(x(\theta)))$. As noted in section 2, the polar wheel is given by $r(\theta) = -y(x(\theta))$ where $x'(\theta) = \cosh x(\theta)$. Separation of variables gives $\theta = 2\tan^{-1}\left(\tanh\frac{x}{2}\right)$ or $x(\theta) = 2\tanh^{-1}\left(\tan\frac{\theta}{2}\right)$. Then $r(\theta) = -\left(-\cosh\left(2\tanh^{-1}\left(\tan\frac{\theta}{2}\right)\right)\right) = \sec\theta$, the polar form of the line $y = -1$.

## 4. Rolling down a hill.

Suppose the road is the tilted line $y = -k\,x - 1$. Then $x'(\theta) = k\,x(\theta) + 1$, giving $x(\theta) = \left(e^{k\theta} - 1\right)/k$ and the wheel is given by $r(\theta) = -(-k\,x(\theta) - 1) = e^{k\theta}$. So a logarithmic spiral rolls smoothly down a straight slope (Figure 9). This is example 2 in the table in [**1**].

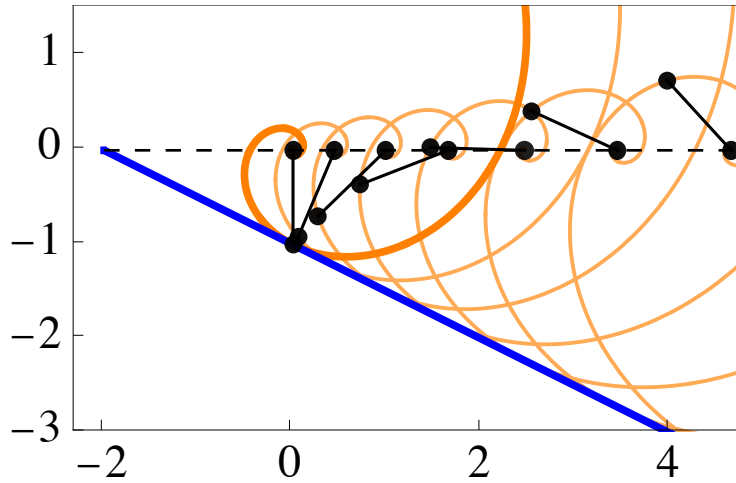

**Figure 9.** The logarithmic spiral $r = e^{k\theta}$ centered at the origin rolls smoothly down the tilted line $y = -k\,x - 1$; here $k = 1/2$.

We can truncate the road and wheel to make a wheel that smoothly rolls along a sawtooth (Figure 10). For four lobes we use the fact that $W\left(\frac{\pi}{4}\right)$ lies on the 45° line to the southwest of the origin. The arc length of the spiral to that point is $\sqrt{2}\,\left(e^{\pi/4} - 1\right)$, and that guides the formation of the teeth forming the road.

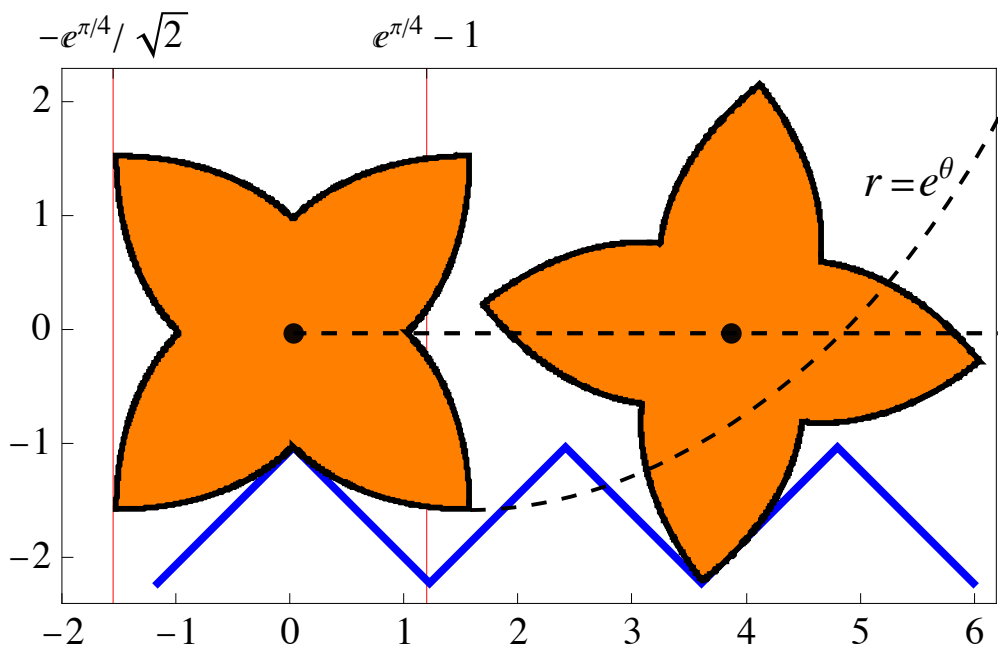

**Figure 10.** A wheel formed from pieces of a logarithmic spiral rolls smoothly on a sawtooth road.

## 5. Road for a nonstandard round wheel.

One can consider a wheel that is a circle of radius 1 whose center is taken to be a point on its circumference. Such is given by $r(\theta) = 2\cos\theta$, with $-\pi/2 < \theta < \pi/2$. The open interval is used here because the theorem applies only when $r(\theta) > 0$. Then $x = \int r(\theta)\, d\theta = 2\sin\theta$, so $\theta = \sin^{-1}(x/2)$ and $y = -2\cos\theta = -\sqrt{4-x^2}$; therefore the road satisfies $x^2 + y^2 = 4$ and is a semicircle of radius 2 (Figure 11). We can extend the rolling by symmetry to the full circular road, and this corresponds to the classic geometrical fact—known as Tusi's mechanism and dating from 1247 CE—that when a wheel rolls without slipping inside a wheel twice as large, any point on the circumference of the rolling wheel traces out a diameter of the larger wheel.

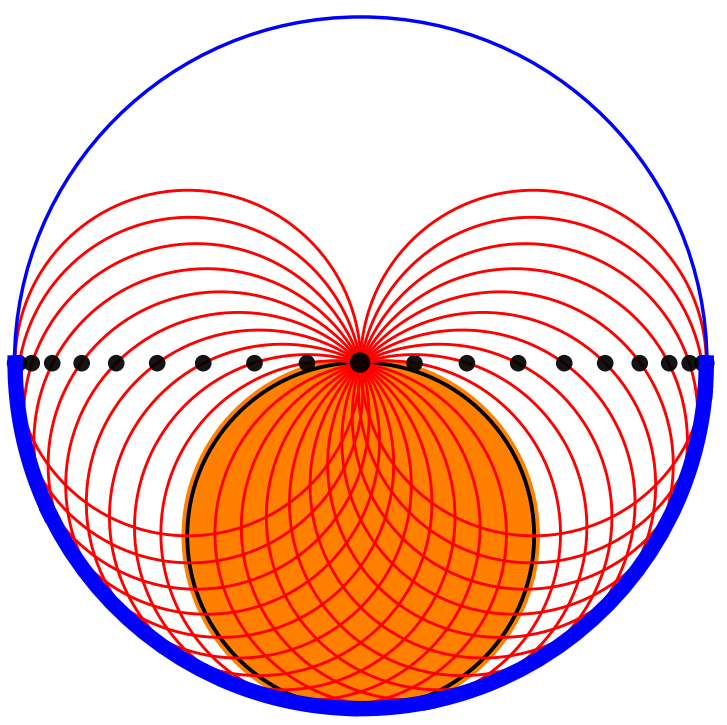

**Figure 11.** A round wheel whose center is taken to be on its circumference rolls smoothly inside a wheel having twice the radius.

## 6. A wheel congruent to its road.

Suppose the wheel is the parabola $y = x^2 - \frac{1}{4}$; in polar form this is $r = \frac{1}{2(1+\cos\theta)}$. The theorem then gives the road as $\left(\frac{1}{2}\tan\frac{\theta}{2}, -\frac{1}{2(1+\cos\theta)}\right)$, which is the same as $y = -x^2 - \frac{1}{4}$ (Figure 12). So we have the remarkable situation of a wheel for which the appropriate road is its reflection; this was discovered by Robison [**1**].

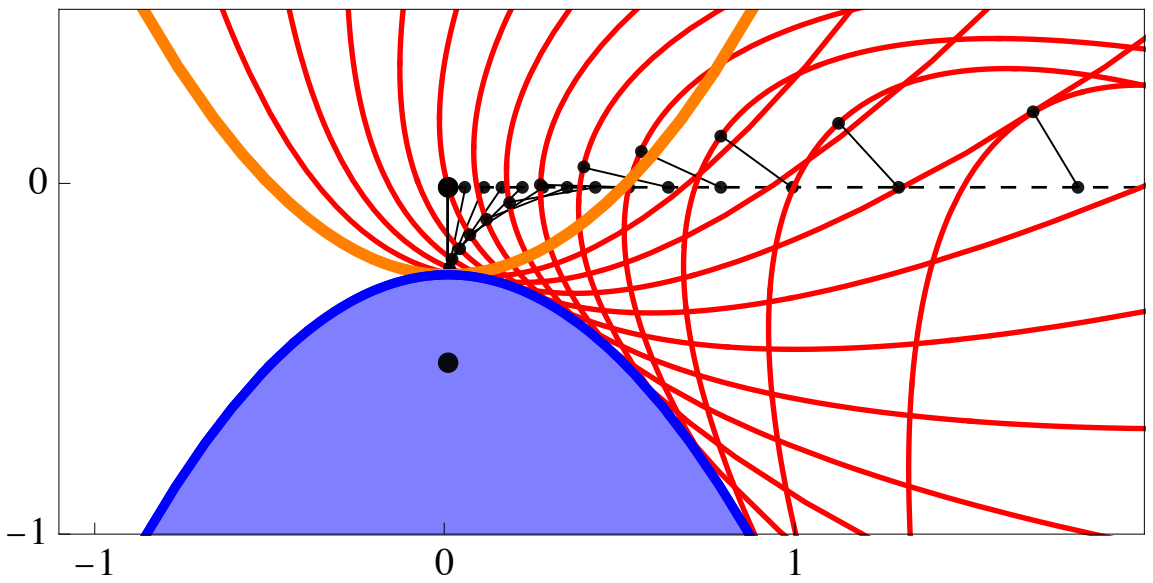

**Figure 12.** If the wheel is a parabola, with its center taken to be its focus, then the road is congruent to the wheel. Here the wheel is $y = x^2 - \frac{1}{4}$.

Robison proved that under certain conditions a wheel that is congruent to its own road must be a parabola. His exact result is as follows.

**Proposition.** If a wheel and a road are congruent and corresponding points under this congruence are in contact when rolling, then the road is a parabola (or part of a parabola) whose directrix is the horizontal axis.

*Proof*. Let $C$ be the congruence that maps the wheel in its initial position onto the road. Let $O = (0, 0)$ and $F = C(O)$. Take an arbitrary point $P'$ on the road, and let $P = C^{-1}(P')$ be the corresponding point on the wheel (Figure 13). By assumption, when the wheel rolls to the contact point $P'$, the point $P$ moves to $P'$. At this moment, the wheel's center is at a point $O'$ on the $x$-axis. Rolling as well as the congruence $C$ preserve lengths, so $|O' P'| = |OP| = \left|C^{-1}(FP')\right| = |FP'|$. By part (b) of the corollary, the segment $O' P'$ is vertical, so its length equals the distance of $P'$ from the $x$-axis. Hence, $P'$ lies on the parabola with focus $F$ and horizontal axis as directrix. □

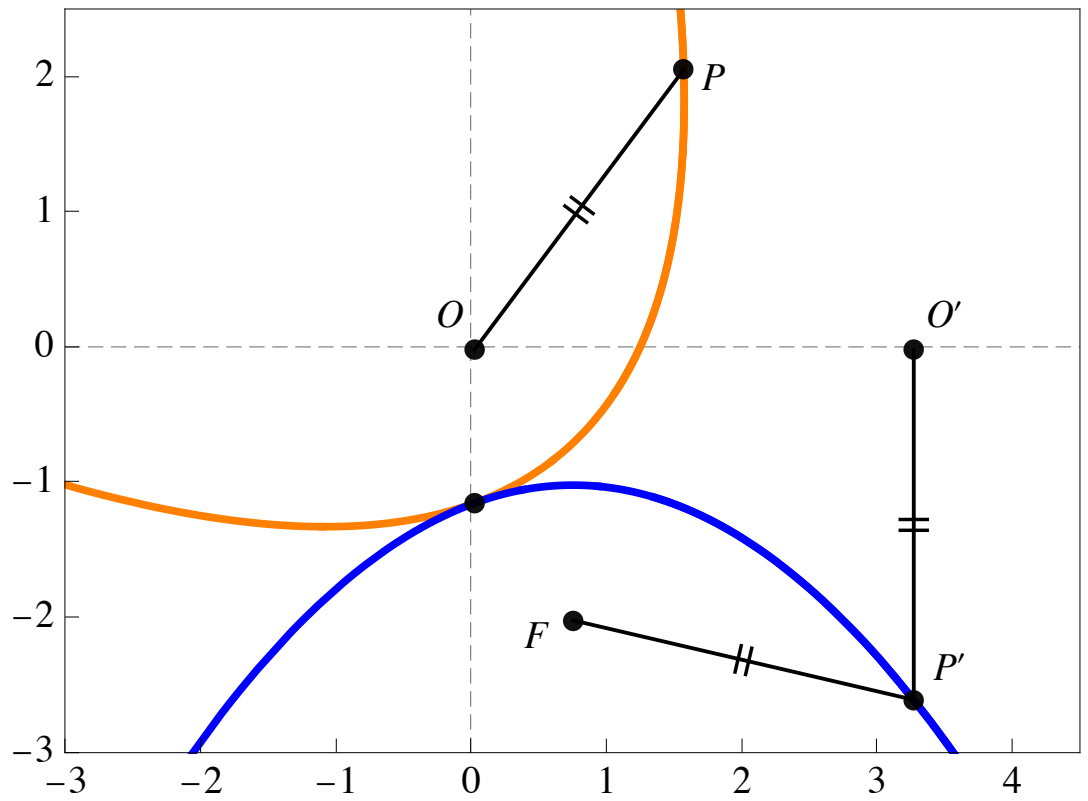

**Figure 13.** The wheel is a rotation of the parabola $y = x^2 - \frac{1}{4}$.

Here are some additional observations:

- If a wheel is congruent to the road as above, it is a parabola. In the initial position its directrix need not be horizontal, but we can always achieve this by rolling the wheel so that its center lies directly above $F$.

- By the previous comment, we can choose the initial position so that $F = C(O)$ lies on the $y$-axis, say at $F = (0, -d)$. Then the vertices of both parabolas touch at $(0, -d/2)$, the equation of the wheel is $y = x^2 / (2d) - d/2$ (or $r = d/(1 + \cos\theta)$ in polar form), and the road is $y = -x^2 / (2d) - d/2$. Up to the choice of the wheel's initial position, these are the only wheel–road pairs that are congruent and such that the corresponding points of this congruence are in contact during rolling.

Robison's condition—that corresponding points are contact points—can be formulated in a different way. We can instead assume the following:

1. The wheel is piecewise continuously differentiable.

2. There exists a point $P$ on the wheel and a point $P'$ on the road such that $P' = C(P)$, and $P$ contacts $P'$ at some stage of rolling.

3. The congruence $C$ is such that increasing the parameter $\theta$ on the wheel increases the $x$-coordinate of the corresponding point on the road.

Under these assumptions, it follows that all points corresponding to each other under $C$ are in contact during rolling. Indeed, since $C$ preserves curve lengths, the image of an arbitrary point $Q$ on the wheel can be determined by measuring its signed distance along the wheel from $P$ and locating a point $Q'$ on the road at the same signed distance from $P'$. Because there is no slipping, $Q$ and its image $Q' = C(Q)$ are contact points.

Note that condition (3) cannot be dropped. For example, if the wheel and road are as in Figure 12), there are two possibilities for $C$, reflection and rotation by 180°. In the latter case, increasing $\theta$ decreases $x$, and the only corresponding points that are also contact points are the vertices of both parabolas. Question: Is it possible to remove condition (2); that is, does the mere fact that the wheel and road are congruent, together with (1), imply the existence of at least one pair of contact points $P$, $P'$ such that $C(P) = P'$? If this could be proved, it would follow that if a road is congruent to its wheel, is piecewise continuously differentiable, and condition (3) holds, then the road is a parabola or part of a parabola.

## 7. A monster wheel.

The theorem requires only that the wheel be continuous; it need not be rectifiable or differentiable. So it is natural to try to roll a wheel that is a Weierstrass monster: a continuous but nowhere differentiable function. The general Weierstrass function is $f(x) = \sum_{n=0}^{\infty} a^{-n} \cos(b^n x)$, where $b > a > 1$. The graph of $f(x)$ is self-similar and is essentially the historically first fractal curve; see Figure 14. The Hausdorff dimension of the graph of $f(x)$ was conjectured to be $2 - (\ln a)/(\ln b)$, a result that was proved relatively recently [**6**]. For a polar wheel we must have positive radius, so we will define a monstrous polar wheel by $r(\theta) = \left(\sum_{n=0}^{\infty} 2^{-n} \cos(3^n \theta)\right) - 3$ (Figure 15).

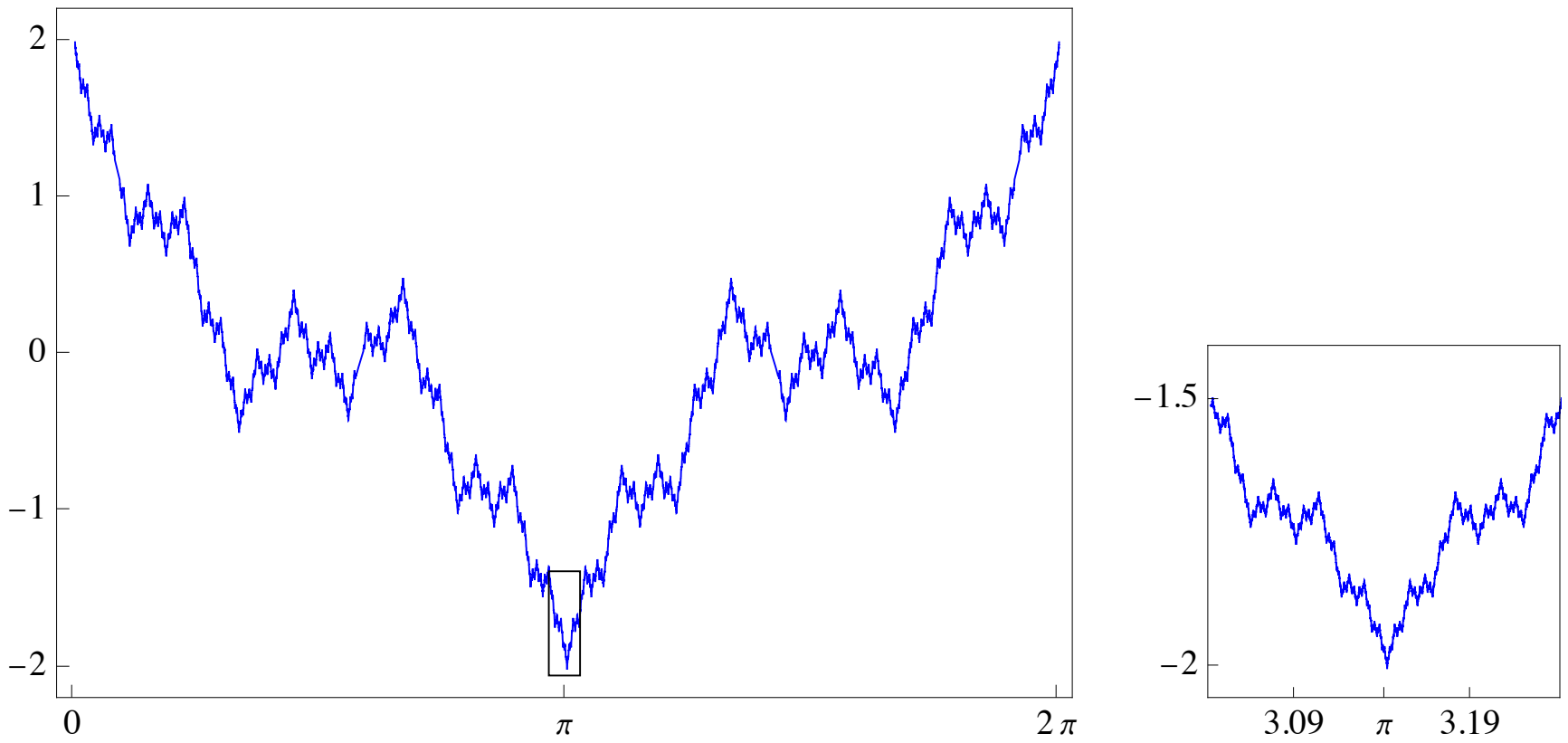

**Figure 14.** The continuous, nowhere differentiable function $\sum_{n=0}^{\infty} 2^{-n} \cos(3^n x)$. The closeup image shows the self-similarity. Its fractal dimension is about 1.37.

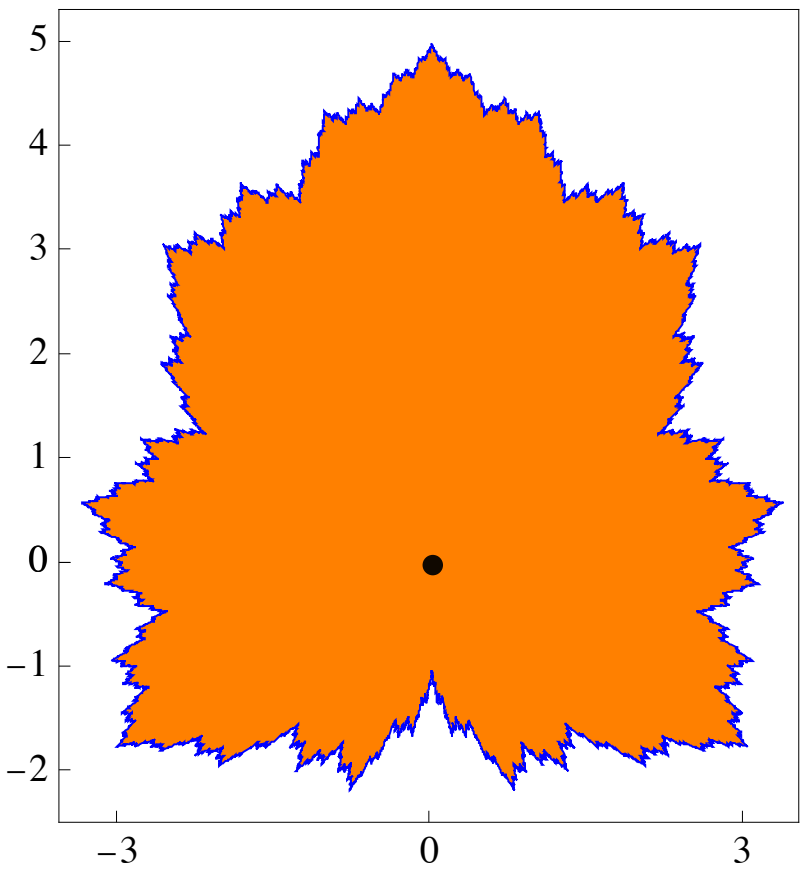

**Figure 15.** A monster wheel.

It is somewhat surprising that the integral of a Weierstrass function is nicely smooth, and that means that the construction of the road as $(x(\theta), y(\theta))$ is straightforward using a numerical approach to the differential equation. Figure 16 shows $x(\theta)$, obtained from a numerical solution to the differential equation $x'(\theta) = r(\theta)$ (using 50 terms in the sum). The graph of $x(\theta)$ looks nice, but its derivative is $r(\theta)$, essentially the same as the Weierstrass function shown in Figure 14. And its second derivative does not exist at any point whatsoever so, for example, it makes no sense to talk about the curvature of $x(\theta)$.

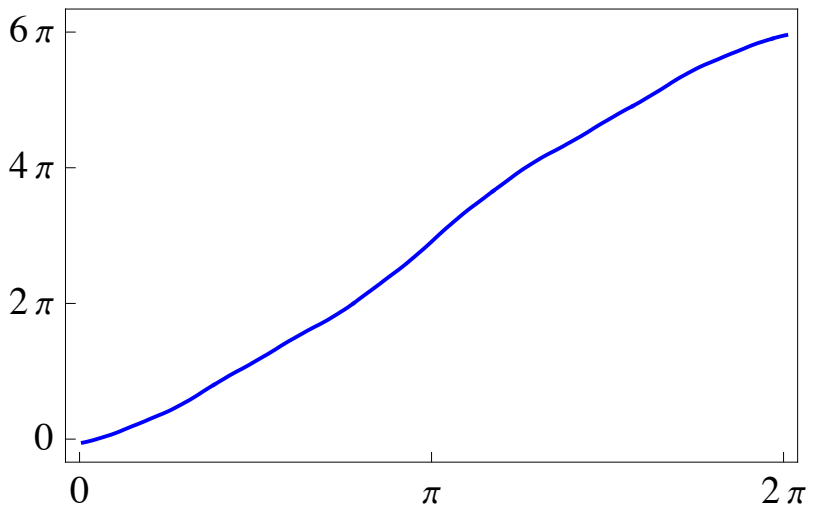

**Figure 16.** The graph of $x(\theta)$ appears remarkably smooth.

We can now construct the road on which the monster wheel will roll. As with a rolling triangle, the wheel will often crash into the road. Indeed, it appears that some part of the wheel is below the road in *every* rolled position. Nevertheless, the road exists and is not hard to compute; Figure 17 shows the rolling wheel and road. The road is not rectifiable, but piecewise linear approximations have lengths that are unbounded, so one can say that both the road and the wheel have infinite length over any nondegenerate interval.

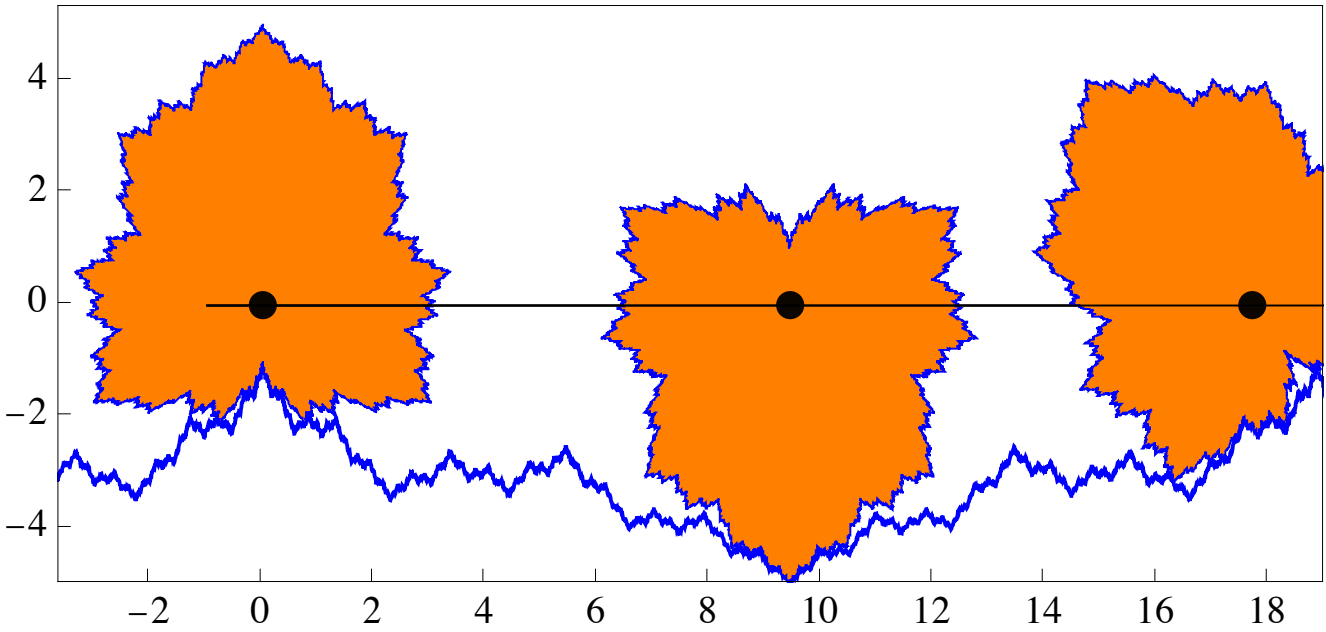

**Figure 17.** A monster wheel rolling along a road that accommodates its nowhere differentiable shape.

The path followed by any point on the wheel is somewhat well behaved; two such paths are shown in Figure 18. The path arising from the initial contact point ($\theta = 0$) never intersects the road, except at its contact point. It appears that all other paths pass below the road at some time during rolling. The smoothness of these paths belies their true nature. Each is differentiable, but the derivatives look like the Weierstrass function. So the loop in the figure is a differentiable curve, with continuous derivative. But there is no second derivative, so the curvature of the curve is not defined anywhere.

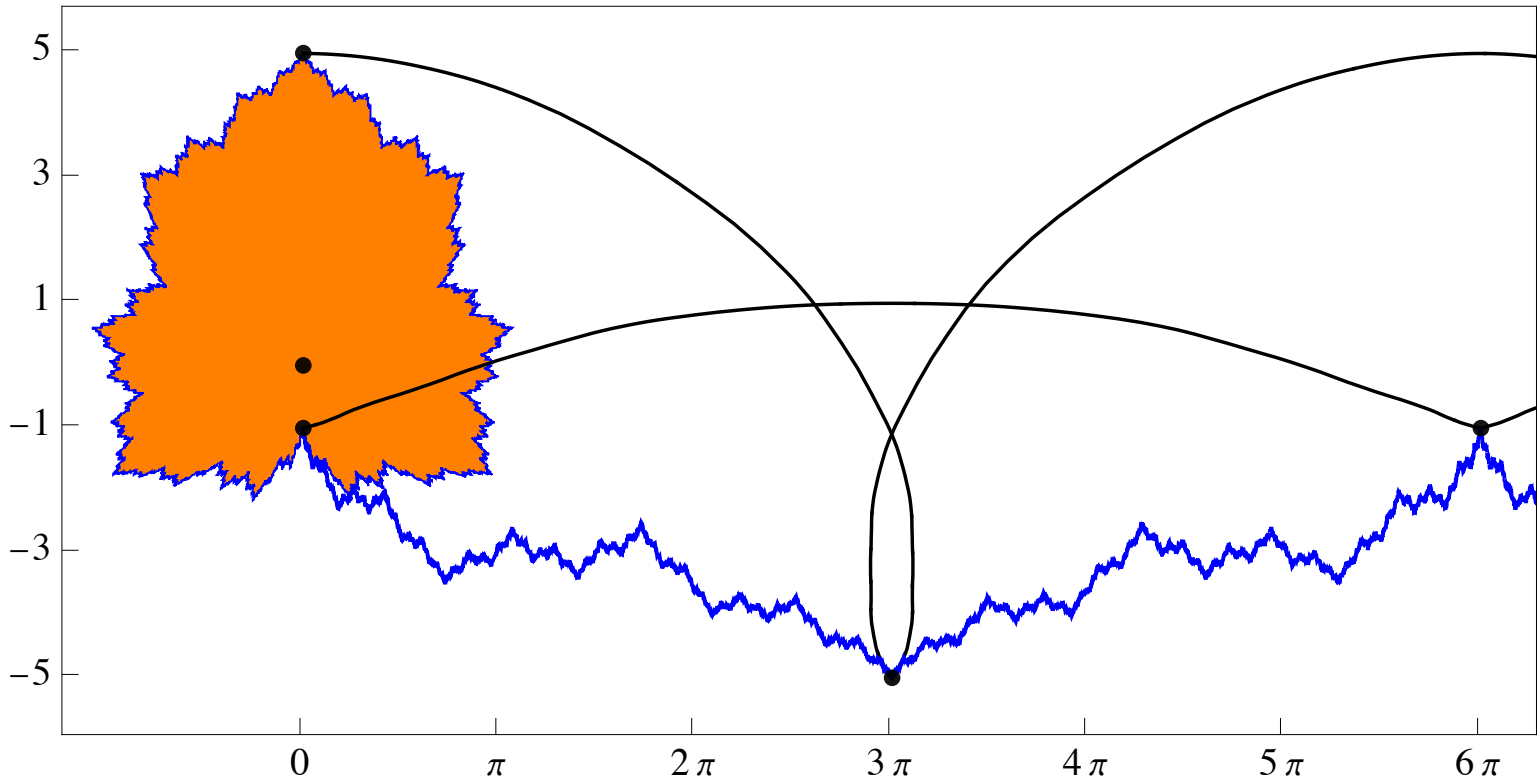

**Figure 18.** The rolling motion for two points of the wheel: $\theta = 0$ and $\theta = \pi$.

## 4. INADEQUACY OF EQUAL ARC LENGTHS.

Corollary (c) states that for piecewise $C^1$ curves, the no-slipping condition implies equality of arc lengths. Given the prominence of the arc-length condition in earlier work, it is natural to ask about the converse; after all, it seems as if arc-length equality would imply no slipping. But it does not. Consider the straight line wheel $y = -1$, with polar form $r = \sec\theta$. If $\rho(\theta) = \theta$, we get the classic catenary road, so try $\rho(\theta) = 2\theta$ (coefficients other than 2 work, but the situation is messier). The road is then $(\tau(2\theta) - \tan\theta, -1)$, which is identical to the wheel, parametrized differently (Figure 19). To satisfy the arc-length condition, suppose the contact point for $W(\theta)$ after some rolling is $C(\theta)$; we seek equality between the length along the wheel to $W(\theta)$, which is $\tan\theta$, the $x$-coordinate of $W(\theta)$, and the length of the road to $C(\theta)$, which is the $x$-coordinate of $M_{2\theta} \cdot W(\theta) + (\tau(2\theta), 0)$, or $\tau(2\theta) - \tan\theta$. Therefore $\tau(2\theta) = 2\tan\theta$, or more generally, $\tau(\phi) = 2\tan(\phi/2)$. And then $C(\theta)$ equals the original point $W(\theta)$. The equality of slopes is clearly false for all $\theta > 0$ and so by Corollary (d) there must be slipping. But the arc-length condition holds. Combining the arc-length condition with another condition does work to define rolling; but our theorem shows that a single no-slipping condition suffices (and rectifiability is not required).

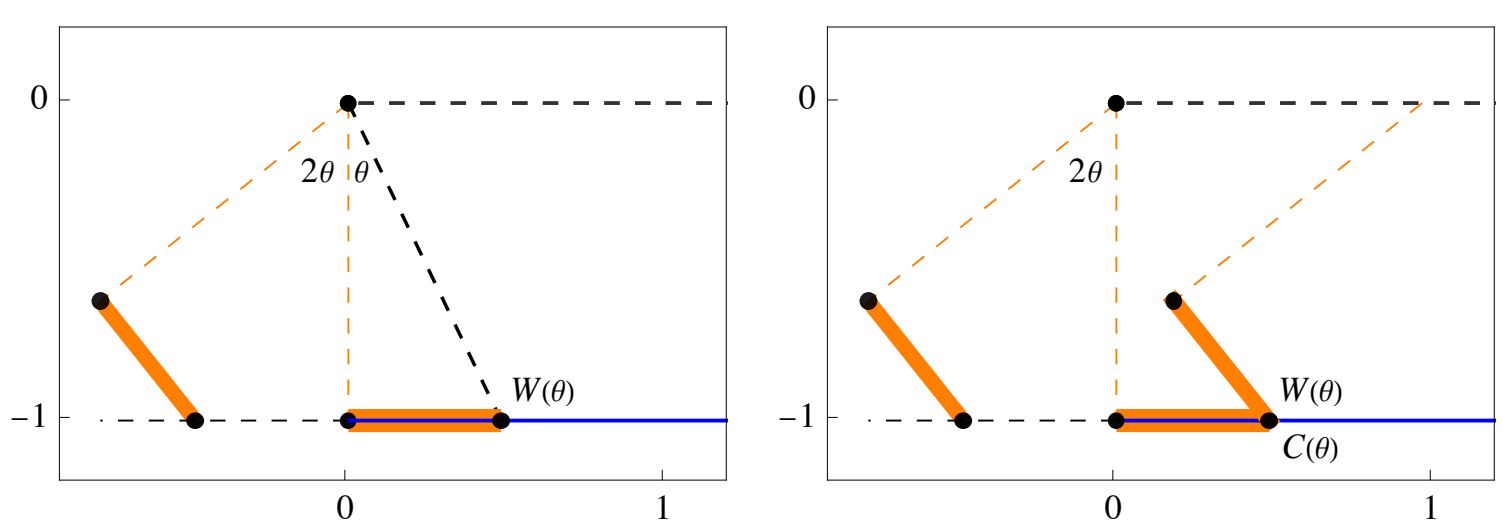

**Figure 19.** A line segment moving along itself. Left: Rotation, but no translation. Right: After translation to the contact point, the relevant arc lengths are equal (isosceles triangle).

## 5. CONCLUSION.

The idea that a square wheel can roll smoothly continues to delight people, from children and adults visiting the National Museum of Mathematics in New York City [**5**], to the bridge designers in England who constructed the whimsical but efficient Cody Dock Bridge that rolls along giant catenaries [**3**]. The development given here derives the key equations from only the minimal assumption that the wheel does not slip while rolling along the road, thus making the underlying mathematical derivation easier to grasp.